\documentclass[12pt,reqno]{amsart}
\usepackage{amssymb,amscd}
\newtheorem{thm}{Theorem}[section]

\newtheorem{lemma}[thm]{Lemma}
\newtheorem{cor}[thm]{Corollary}
\newtheorem{definition}[thm]{Definition}

\def\bN{{\mathbb{N}}}
\def\bR{{\mathbb{R}}}
\def\bP{{\mathbb{P}}}
\def\bE{{\mathbb{E}}}
\def\sym{{\mathrm{sym}}}
\def\cX{{\mathcal{X}}}
\def\cY{{\mathcal{Y}}}
\def\eps{\varepsilon}
\def\bx{{\mathbf{x}}}
\def\by{{\mathbf{y}}}
\def\bs{{\mathbf{s}}}
\def\cT{{\mathcal{T}}}

\begin{document}
\title[A new approach to mutual information]
{A new approach to mutual information}

\author[F. Hiai]{Fumio Hiai$\,^{1}$}
\address{Graduate School of Information Sciences,
Tohoku University, Aoba-ku, Sendai 980-8579, Japan}
\author[D. Petz]{D\'enes Petz$\,^{2}$}
\address{Alfr\'ed R\'enyi Institute of Mathematics, Hungarian 
Academy of Sciences, H-1053 Budapest, Re\'altanoda u. 13-15, Hungary}

\thanks{$^1\,$Supported in part by Grant-in-Aid for Scientific Research
(B)17340043.}
\thanks{$^2\,$Supported in part by the Hungarian Research Grant OTKA T068258.}


\maketitle

\begin{abstract}
A new expression as a certain asymptotic limit via ``discrete micro-states" of
permutations is provided to the mutual information of both continuous and discrete random variables.
\end{abstract}

\section*{Introduction}

One of the important quantities in information theory is the mutual 
information of two random variables $X$ and $Y$ which is expressed
in terms of the Boltzmann-Gibbs entropy $H(\cdot)$ as follows:
$$
I(X\wedge Y)=-H(X,Y)+H(X)+H(Y)
$$
when $X,Y$ are continuous variables. For the expression of $I(X\wedge Y)$ of discrete
variables $X,Y$, the above $H(\cdot)$ is replaced by the Shannon entropy. A more
practical and rigorous definition via the relative entropy is
$$
I(X \wedge Y):=S(\mu_{(X,Y)},\mu_X\otimes\mu_Y),
$$
where $\mu_{(X,Y)}$ denotes the joint distribution measure of $(X,Y)$ and
$\mu_X\otimes\mu_Y$ the product of the respective distribution measures 
of $X,Y$.

The aim of this paper is to show that the mutual information $I(X\wedge Y)$ 
is gained as a certain asymptotic limit of the volume of ``discrete 
micro-states" consisting of permutations approximating joint moments of 
$(X,Y)$ in some way. In Section 1, more generally we consider an $n$-tuple of
real  bounded random variables $(X_1,\dots,X_n)$. Denote by
$\Delta(X_1,\dots,X_n;N,m,\delta)$ the set of $(\bx_1,\dots,\bx_n)$ of 
$\bx_i\in\bR^N$ whose joint moments (on the uniform distributed $N$-point set) 
of order up to $m$ approximate those of $(X_1,\dots,X_n)$ up to an error 
$\delta$. Furthermore, denote by $\Delta_\sym(X_1,\dots,X_n;N,m,\delta)$ 
the set of $(\sigma_1,\dots,\sigma_n)$ of permutations $\sigma_i\in S_N$ 
such that
$(\sigma_1(\bx_1),\dots,\sigma_n(\bx_n))\in\Delta(X_1,\dots,X_n;N,m,\delta)$
for some $\bx_1,\dots,\bx_n\in\bR_\le^N$, where $\bR_\le^N$ is the $\bR^N$-vectors
arranged in increasing order. Then, the asymptotic volume
$$
{1\over N}\log
\gamma_{S_N}^{\otimes n}\bigl(\Delta_\sym(X_1,\dots,X_n;N,m,\delta)\bigr)
$$
under the uniform probability measure $\gamma_{S_N}$ on $S_N$ is shown to 
converge as $\limsup_{N\to\infty}$ (also $\liminf_{N\to\infty}$) and then
$\lim_{m\to\infty,\delta\searrow0}$ to
$$
-H(X_1,\dots,X_n)+\sum_{i=1}^nH(X_i)
$$
as long as $H(X_i)>-\infty$ for $1\le i\le n$. Thus, we obtain a kind of
discretization of the mutual information via symmetric group (or permutations).

The approach can be applied to an $n$-tuple of discrete random variables
$(X_1,\dots,X_n)$ as well. But the definition of the $\Delta_\sym$-set of
micro-states for discrete variables is somewhat different from the continuous
variable case mentioned above, and we discuss the discrete variable case in
Section 2 separately.

The idea comes from the paper \cite{HMU}. Motivated by theory of mutual free
information in \cite{V6}, a similar approach to Voiculescu's free entropy
is provided there. The free entropy is the free probability counterpart of
the Boltzmann-Gibbs entropy, and $\bR^N$-vectors and the symmetric group
$S_N$ here are replaced by Hermitian $N\times N$ matrices and the unitary
group $\mathrm{U}(N)$, respectively. In this way, the ``discretization
approach" here is in some sense a classical analog of the ``orbital approach"
in \cite{HMU}.

\section{The continuous case}

For $N\in\bN$ let $\bR_\le^N$ be the convex cone of the $N$-dimensional Euclidean
space $\bR^N$ consisting of $\bx=(x_1,\dots,x_N)$ such that
$x_1\le x_2\le\dots\le x_N$. The space $\bR^N$ is naturally regarded as the real
function algebra on the $N$-point set. Let $S_N$ be the symmetric group of order
$N$ (i.e., the permutations on $\{1,2,\dots,n\}$). Throughout this section let
$(X_1,\dots,X_n)$ be an $n$-tuple of real random variables on a probability space
$(\Omega,\bP)$, and assume that the $X_i$'s are bounded (i.e.,
$X_i\in L^\infty(\Omega;\bP)$). The {\it Boltzmann-Gibbs entropy} of
$(X_1,\dots,X_n)$ is defined to be
$$
H(X_1,\dots,X_n):=-\int\cdots\int_{\bR^n}p(x_1,\dots,x_n)\log p(x_1,\dots,x_n)
\,dx_1\cdots dx_n
$$
if the joint density $p(x_1,\dots,x_n)$ of $(X_1,\dots,X_n)$ exists; otherwise
$H(X_1,\dots,X_n)=-\infty$. Note that the above integral is well defined in
$[-\infty,\infty)$ since the density $p$ is compactly supported.

\begin{definition}\label{D-1.1}{\rm
The mean value of $\bx=(x_1,\dots,x_N)$ in $\bR^N$ is given by
$$
\kappa_N(\bx):={1\over N}\sum_{j=1}^Nx_j.
$$
For each $N,m\in\bN$ and $\delta>0$ we define $\Delta(X_1,\dots,X_n;N,m,\delta)$ to
be the set of all $n$-tuples $(\bx_1,\dots,\bx_n)$ of
$\bx_i=(x_{i1},\dots,x_{iN})\in\bR^N$, $1\le i\le n$, such that
$$
|\kappa_N(\bx_{i_1}\cdots\bx_{i_k})-\bE(X_{i_1}\cdots X_{i_k})|<\delta
$$
for all $1\le i_1,\dots,i_k\le n$ with $1\le k\le m$, where
$\bx_{i_1}\cdots\bx_{i_k}$ means the pointwise product, i.e.,
$$
\bx_{i_1}\cdots\bx_{i_k}:=(x_{i_11}\cdots x_{i_k1},
x_{i_12}\cdots x_{i_k2},\dots,x_{i_1N}\cdots x_{i_kN})\in\bR^N
$$
and $\bE(\cdot)$ denotes the expectation on $(\Omega,\bP)$. For each $R>0$, define
$\Delta_R(X_1,\dots,X_n;\allowbreak N,m,\delta)$ to be the set of all
$(\bx_1,\dots,\bx_n)\in\Delta(X_1,\dots,X_n;N,m,\delta)$ such that
$\bx_i\in[-R,R]^N$ for all $1\le i\le n$.
}\end{definition}

Heuristically, $\Delta(X_1,\dots,X_n;N,m,\delta)$ is the set of ``micro-states"
consisting of $n$-tuples of discrete random variables on the $N$-point set with
the uniform probability such that all joint moments of order up to $m$ give the
corresponding joint moments of $X_1,\dots,X_n$ up to an error $\delta$.

For $\bx\in\bR^N$ write $\|\bx\|_p:=(N^{-1}\sum_{j=1}^N|x_j|^p)^{1/p}$ for
$1\le p<\infty$ and $\|\bx\|_\infty:=\max_{1\le j\le N}|x_j|$ while $\|X\|_p$ denotes
the $L^p$-norm of a real random variable $X$ on $(\Omega,\bP)$.

The next lemma is seen from \cite[5.1.1]{HP} based on the Sanov large deviation
theorem, which says that the Boltzmann-Gibbs entropy is gained as an asymptotic limit
of the volume of the approximating micro-states.

\begin{lemma}\label{L-1.2}
For every $m\in\bN$ and $\delta>0$ and for any choice of
$R\ge\max_{1\le i\le n}\|X_i\|_\infty$, the limit
$$
\lim_{N\to\infty}{1\over N}\log
\lambda_N^{\otimes n}\bigl(\Delta_R(X_1,\dots,X_n;N,m,\delta)\bigr)
$$
exists, where $\lambda_N$ is the Lebesgue measure on $\bR^N$. Furthermore, one has
$$
H(X_1,\dots,X_n)=\lim_{m\rightarrow\infty,\delta\searrow 0}
\lim_{N\to\infty}{1\over N}\log
\lambda_N^{\otimes n}\bigl(\Delta_R(X_1,\dots,X_n;N,m,\delta)\bigr)
$$
independently of the choice of $R\ge\max_{1\le i\le n}\|X_i\|_\infty$.
\end{lemma}

In the following let us introduce some kinds of mutual information in the
discretization approach using micro-states of permutations.

\begin{definition}\label{D-1.3}{\rm
The action of $S_N$ on $\bR^N$ is given by
$$
\sigma(\bx):=(x_{\sigma^{-1}(1)},x_{\sigma^{-1}(2)},\dots,x_{\sigma^{-1}(N)})
$$
for $\sigma\in S_N$ and $\bx=(x_1,\dots,x_N)\in\bR^N$. For each $N,m\in\bN$,
$\delta>0$ and $R>0$ we denote by $\Delta_{\sym,R}(X_1,\dots,X_n;N,m,\delta)$ the
set of all $(\sigma_1,\dots,\sigma_n)\in S_N^n$ such that
$$
(\sigma_1(x_1),\dots,\sigma_n(x_n))\in\Delta_R(X_1,\dots,X_n;N,m,\delta)
$$
for some $(\bx_1,\dots,\bx_n)\in(\bR_\le^N)^n$. For each $R>0$ define
$$
I_{\sym,R}(X_1,\dots,X_n):=-\lim_{m\to\infty,\delta\searrow0}
\limsup_{N\to\infty}{1\over N}\log
\gamma_{S_N}^{\otimes n}\bigl(\Delta_{\sym,R}(X_1,\dots,X_n;N,m,\delta)\bigr),
$$
where $\gamma_{S_N}$ is the uniform probability measure on $S_N$. Define also
$\overline{I}_{\sym,R}(X_1,\dots,X_n)$ by replacing $\limsup$ by $\liminf$. Obviously,
$$
0\le I_{\sym,R}(X_1,\dots,X_n)\le\overline{I}_{\sym,R}(X_1,\dots,X_n).
$$
Moreover, $\Delta_{\sym,\infty}(X_1,\dots,X_n;N,m,\delta)$ is defined by replacing
$\Delta_R(X_1,\dots,X_n;N,m,\delta)$ in the above by
$\Delta(X_1,\dots,X_n;N,m,\delta)$ without cut-off by the parameter $R$. Then
$I_{\sym,\infty}(X_1,\dots,X_n)$ and $\overline{I}_{\sym,\infty}(X_1,\dots,X_n)$ are
also defined as above.
}\end{definition}

\begin{definition}\label{D-1.4}{\rm
For each $1\le i\le n$ we choose and fix a sequence $\xi_i=\{\xi_i(N)\}$ of
$\xi_i(N)\in\bR_\le^N$, $N\in\bN$, such that $\kappa_N(\xi_i(N)^k)\to\bE(X_i^k)$ as
$N\to\infty$ for all $k\in\bN$, i.e., $\xi_i(N)\to X_i$ in moments. For each
$N,m\in\bN$ and $\delta>0$ we define
$\Delta_\sym(X_1,\dots,X_n:\xi_1(N),\dots,\xi_n(N);N,m,\delta)$ to be the set of all
$(\sigma_1,\dots,\sigma_n)\in S_N^n$ such that
$$
(\sigma_1(\xi_1(N)),\dots,\sigma_n(\xi_n(N)))\in\Delta(X_1,\dots,X_n;N,m,\delta).
$$
Define
\begin{align*}
&I_\sym(X_1,\dots,X_n:\xi_1,\dots,\xi_n) \\
&\quad:=-\lim_{m\rightarrow\infty,\delta\searrow 0}
\limsup_{N\to\infty}{1\over N}\log\gamma_{S_N}^{\otimes n}
\bigl(\Delta_\sym(X_1,\dots,X_n:\xi_1(N),\dots,\xi_n(N);N,m,\delta)\bigr)
\end{align*}
and $\overline{I}_\sym(X_1,\dots,X_n:\xi_1,\dots\xi_n)$ by replacing $\limsup$ by
$\liminf$.
}\end{definition}

The next proposition asserts that the quantities in Definitions \ref{D-1.3} and
\ref{D-1.4} are all equivalent.

\begin{lemma}\label{L-1.5}
For any choice of $R\ge\max_{1\le i\le n}\|X_i\|_\infty$ and for any choices of
approximating sequences $\xi_1,\dots,\xi_n$ one has
\begin{align}
I_{\sym,\infty}(X_1,\dots,X_n)&=I_{\sym,R}(X_1,\dots,X_n)
=I_\sym(X_1,\dots,X_n:\xi_1,\dots,\xi_n), \label{F-1.1}\\
\overline{I}_{\sym,\infty}(X_1,\dots,X_n)&=\overline{I}_{\sym,R}(X_1,\dots,X_n)
=\overline{I}_\sym(X_1,\dots,X_n:\xi_1,\dots,\xi_n). \label{F-1.2}
\end{align}
\end{lemma}

\begin{proof}
It is obvious that $\Delta_\sym(X_1,\dots,X_n:\xi_1(N),\dots,\xi_n(N);N,m,\delta)$
is included in $\Delta_{\sym,\infty}(X_1,\dots,X_n;N,m,\delta)$ for any approximating
sequences $\xi_i$. Moreover, for each $1\le i\le n$ an approximating sequence $\xi_i$
can be chosen so that $\|\xi_i(N)\|_\infty\le\|X_i\|_\infty$ for all $N$; then
$\Delta_\sym(X_1,\dots,X_n:\xi_1(N),\dots,\xi_n(N);N,m,\delta)\subset
\Delta_{\sym,R}(X_1,\dots,X_n;\allowbreak N,m,\delta)$ for any
$R\ge R_0:=\max_{1\le i\le n}\|X_i\|_\infty$. Hence it suffices to prove that for any
approximating sequences $\xi_i$ and for every $m\in\bN$ and $\delta>0$, there are an
$m'\in\bN$, a $\delta'>0$ and an $N_0\in\bN$ so that
$$
\Delta_{\sym,\infty}(X_1,\dots,X_n;N,m',\delta')
\subset\Delta_\sym(X_1,\dots,X_n:\xi_1(N),\dots,\xi_n(N);N,m,\delta)
$$
for all $N\ge N_0$. Choose a $\rho\in(0,1)$ with $m(R_0+1)^{m-1}\rho<\delta/2$.
By \cite[Lemma 4.3]{V2} (also \cite[4.3.4]{HP}) there exist an $m'\in\bN$ with
$m'\ge2m$, a $\delta'>0$ with $\delta'\le\min\{1,\delta/2\}$ and an $N_0\in\bN$
such that for every $1\le i\le n$ and every $\bx\in\bR_\le^N$ with $N\ge N_0$, if
$|\kappa_N(\bx^k)-\bE(X_i^k)|<\delta'$ for all $1\le k\le m'$, then
$\|\bx-\xi_i(N)\|_m<\rho$. Suppose $N\ge N_0$ and
$(\sigma_1,\dots,\sigma_n)\in\Delta_{\sym,\infty}(X_1,\dots,X_n;N,m',\delta')$;
then $(\sigma_1(\bx_1),\dots,\sigma_n(\bx_n))\in\Delta(X_1,\dots,X_n;N,m',\delta')$
for some $(\bx_1,\dots,\bx_n)\in(\bR_\le^N)^n$. Since
$|\kappa_N(\bx_i^k)-\bE(X_i^k)|<\delta'$ for all $1\le k\le m'$, we get
$\|\bx_i-\xi_i(N)\|_m\le\rho$ and
\begin{align*}
\|\bx_i\|_m&\le\|\bx_i\|_{2m}=\kappa_N(\bx_i^{2m})^{1/2m} \\
&<(\bE(X_i^{2m})+1)^{1/2m} \\
&\le(R_0^{2m}+1)^{1/2m}\le R_0+1.
\end{align*}
Therefore,
\begin{align*}
&|\kappa_N(\sigma_{i_1}(\xi_{i_1}(N))\cdots\sigma_{i_k}(\xi_{i_k}(N)))
-\bE(X_{i_1}\cdots X_{i_k})| \\
&\quad\le|\kappa_N(\sigma_{i_1}(\xi_{i_1}(N))\cdots\sigma_{i_k}(\xi_{i_k}(N)))
-\kappa_N(\sigma_{i_1}(\bx_{i_1})\cdots\sigma_{i_k}(\bx_{i_k}))| \\
&\qquad+|\kappa_N(\sigma_{i_1}(\bx_{i_1})\cdots\sigma_{i_k}(\bx_{i_k}))
-\bE(X_{i_1}\cdots X_{i_k})| \\
&\quad\le m(R_0+1)^{m-1}\rho+\delta'<\delta
\end{align*}
for all $1\le i_1,\dots,i_k\le n$ with $1\le k\le m$. The above latter inequality
follows from the H\"older inequality. Hence $(\sigma_1,\dots,\sigma_n)\in
\Delta_\sym(X_1,\dots,X_n:\xi_1(N),\dots,\xi_n(N);N,m,\delta)$, and the result
follows.
\end{proof}

Consequently, we denote all the quantities in \eqref{F-1.1} by the same
$I_\sym(X_1,\dots,X_n)$ and those in \eqref{F-1.2} by $\overline{I}_\sym(X_1,\dots,X_n)$.
We call $I_\sym(X_1,\dots,X_n)$ and $\overline{I}_\sym(X_1,\dots,X_n)$ the
{\it mutual information} and {\it upper mutual information} of $(X_1,\dots,X_n)$,
respectively. The terminology ``mutual information" will be justified after the next
theorem.

In the continuous variable case, our main result is the following exact relation of
$I_\sym$ and $\overline{I}_\sym$ with the Boltzmann-Gibbs entropy $H(\cdot)$, which
says that $I_\sym(X_1,\dots,X_n)$ is formally the sum of the separate entropies
$H(X_i)$'s minus the compound $H(X_1,\dots,X_n)$. Thus, a naive meaning of
$I_\sym(X_1,\dots,X_n)$ is the entropy (or information) overlapping among the $X_i$'s.

\begin{thm}\label{T-1.6}
\begin{align*}
H(X_1,\dots,X_n)&=-I_\sym(X_1,\dots,X_n)+\sum_{i=1}^nH(X_i) \\
&=-\overline{I}_\sym(X_1,\dots,X_n)+\sum_{i=1}^nH(X_i).
\end{align*}
\end{thm}

\begin{proof}
If the coordinates $s_i$ of $\bs\in\bR^N$ are all distinct, then $\bs$ is uniquely
written as $\bs=\sigma(\bx)$ with $\bx\in\bR_\le^N$ and $\sigma\in S_N$. Note that
the set of $\bs\in\bR^N$ with $s_i=s_j$ for some $i\ne j$ is a closed subset of
$\lambda_N$-measure zero. Under the correspondence
$$
\bs\in\bR^N\longleftrightarrow(\bx,\sigma)\in\bR_\le^N\times S_N,
\quad\bs=\sigma(\bx)
$$
(well defined on a co-negligible subset of $\bR^N$), the measure $\lambda_N$ is
transformed into the product of $\lambda_N|_{\bR_\le^N}$ and the counting measure on
$S_N$.

In the following proof we adopt, due to Lemma \ref{L-1.5}, the description of
$I_\sym$ and $\overline{I}_\sym$ as $I_{\sym,R}(X_1,\dots,X_n)$ and
$\overline{I}_{\sym,R}(X_1,\dots,X_n)$ with $R:=\max_{1\le i\le n}\|X_i\|_\infty$.
For each $N,m\in\bN$ and $\delta>0$, suppose
$(\bs_1,\dots,\bs_n)\in\Delta_R(X_1,\dots,X_n;N,m,\delta)$ and write
$\bs_i=\sigma_i(\bx_i)$ with $\bx_i\in\bR_\le^N$ and $\sigma_i\in S_N$. Then it is
obvious that
\begin{align*}
&(\bx_1,\dots,\bx_n;\sigma_1,\dots,\sigma_n) \\
&\qquad\in\Biggl(\prod_{i=1}^n\bigl(\Delta_R(X_i;N,m,\delta)\cap\bR_\le^N\bigr)
\Biggr)\times\Delta_{\sym,R}(X_1,\dots,X_n;N,m,\delta).
\end{align*}
By Lemma \ref{L-1.2} and the fact stated at the beginning of the proof, we obtain
\allowdisplaybreaks{
\begin{align*}
H(X_1,\dots,X_n)
&\le\lim_{N\to\infty}{1\over N}
\log\lambda_N^{\otimes n}\bigl(\Delta_R(X_1,\dots,X_n;N,m,\delta)\bigr) \\
&\le\liminf_{N\to\infty}{1\over N}\Biggl(\sum_{i=1}^n
\log\lambda_N\bigl(\Delta_R(X_i;N,m,\delta)\cap\bR_\le^N\bigr) \\
&\qquad\qquad\qquad
+\log\#\Delta_{\sym,R}(X_1,\dots,X_n;N,m,\delta)\Biggr) \\
&=\liminf_{N\to\infty}{1\over N}\Biggl(\sum_{i=1}^n
\log\lambda_N\bigl(\Delta_R(X_i;N,m,\delta)\bigr)-n\log N! \\
&\qquad\qquad\qquad
+\log\#\Delta_{\sym,R}(X_1,\dots,X_n;N,m,\delta)\Biggr) \\
&=\sum_{i=1}^n\lim_{N\to\infty}{1\over N}\log
\lambda_N\bigl(\Delta_R(X_i;N,m,\delta)\bigr) \\
&\qquad
+\liminf_{N\to\infty}{1\over N}
\log\gamma_{S_N}^{\otimes n}\bigl(\Delta_{\sym,R}(X_1,\dots,X_n;N,m,\delta)\bigr).
\end{align*}
}This implies that
\begin{equation}\label{F-1.3}
H(X_1,\dots,X_n)\le\sum_{i=1}^nH(X_i)-\overline{I}_\sym(X_1,\dots,X_n).
\end{equation}

Conversely, for each $m\in\bN$ and $\delta>0$, by \cite[Lemma 4.3]{V2} (also
\cite[4.3.4]{HP}) there are an $m'\in\bN$ with $m'\ge m$, a $\delta'>0$ with
$\delta'\le\delta/2$ and an $N_0\in\bN$ such that for every $N\in\bN$ and
for every $\bx,\by\in\bR_\le^N$, if $\|\bx\|_\infty\le R$ and
$|\kappa_N(\bx^k)-\kappa_N(\by^k)|<2\delta'$ for all $1\le k\le m'$, then
$\|\bx-\by\|_1<\delta/2m(R+1)^{m-1}$. Suppose $N\ge N_0$ and
\begin{align*}
&(\bx_1,\dots,\bx_n;\sigma_1,\dots,\sigma_n) \\
&\qquad\in\Biggl(\prod_{i=1}^n\bigl(\Delta_R(X_i;N,m',\delta')\cap\bR_\le^N\bigr)
\Biggr)\times\Delta_{\sym,R}(X_1,\dots,X_n;N,m',\delta')
\end{align*}
so that
$(\sigma_1(\by_1),\dots,\sigma_n(\by_n))\in\Delta_R(X_1,\dots,X_n;N,m',\delta')$
for some $(\by_1,\dots,\by_n)\in(\bR_\le^N)^n$. Since
$$
|\kappa_N(\bx_i^k)-\kappa_N(\by_i^k)|
\le|\kappa_N(\bx_i^k)-\bE(X_i^k)|+|\kappa_N(\by_i^k)-\bE(X_i^k)|<2\delta'
$$
for all $1\le k\le m'$, we get $\|\bx_i-\by_i\|_1<\delta/2m(R+1)^{m-1}$ for
$1\le i\le n$. Therefore,
\begin{align*}
&|\kappa_N(\sigma_{i_1}(\bx_{i_1})\cdots\sigma_{i_k}(\bx_{i_k}))
-\bE(X_{i_1}\cdots X_{i_k})| \\
&\qquad\le|\kappa_N(\sigma_{i_1}(\bx_{i_1})\cdots\sigma_{i_k}(\bx_{i_k}))
-\kappa_N(\sigma_{i_1}(\by_{i_1})\cdots\sigma_{i_k}(\by_{i_k}))| \\
&\qquad\quad+|\kappa_N(\sigma_{i_1}(\by_{i_1})\cdots\sigma_{i_k}(\by_{i_k}))
-\bE(X_{i_1}\cdots X_{i_k})| \\
&\qquad\le m(R+1)^{m-1}\max_{1\le i\le n}\|\bx_i-\by_i\|_1+\delta' \\
&\qquad<{\delta\over2}+\delta'\le\delta
\end{align*}
for all $1\le i_1,\dots, i_k\le n$ with $1\le k\le m$. This implies that
$(\sigma_1(\bx_1),\dots,\sigma_n(\bx_n))\in\Delta_R(X_1,\dots,X_n;N,m,\delta)$. By
Lemma \ref{L-1.2} we obtain
\allowdisplaybreaks{
\begin{align*}
&\sum_{i=1}^nH(X_i)-I_\sym(X_1,\dots,X_n) \\
&\qquad\le\sum_{i=1}^n\lim_{N\to\infty}{1\over N}
\log\lambda_N\bigl(\Delta_R(X_i;N,m',\delta')\bigr) \\
&\qquad\quad+\limsup_{N\to\infty}{1\over N}\log\gamma_{S_N}^{\otimes n}
\bigl(\Delta_{\sym,R}(X_1,\dots,X_n;N,m',\delta')\bigr) \\
&\qquad=\limsup_{N\to\infty}{1\over N}\Biggl(\sum_{i=1}^n
\log\lambda_N\bigl(\Delta_R(X_i;N,m',\delta')\cap\bR_\le^N\bigr) \\
&\qquad\qquad\qquad\qquad
+\log\#\Delta_{\sym,R}(X_1,\dots,X_n;N,m',\delta')\Biggr) \\
&\qquad\le\limsup_{N\to\infty}{1\over N}
\log\lambda_N^{\otimes n}\bigl(\Delta_R(X_1,\dots,X_n;N,m,\delta)\bigr).
\end{align*}
}This implies by Lemma \ref{L-1.2} once again that
\begin{equation}\label{F-1.4}
\sum_{i=1}^nH(X_i)-I_\sym(X_1,\dots,X_n)\le H(X_1,\dots,X_n).
\end{equation}
The result follows from \eqref{F-1.3} and \eqref{F-1.4}.
\end{proof}

Let $\mu_{(X_1,\dots,X_n)}$ be the joint distribution measure on $\bR^n$ of
$(X_1,\dots,X_n)$ while $\mu_{X_i}$ is that of $X_i$ for $1\le i\le n$. Let
$S(\mu_{(X_1,\dots,X_n)},\mu_{X_1}\otimes\dots\otimes\mu_{X_n})$ denote the
{\it relative entropy} (or {\it the Kullback-Leibler divergence}) of
$\mu_{(X_1,\dots,X_n)}$ with respect to the product measure
$\mu_{X_1}\otimes\dots\otimes\mu_{X_n}$, i.e.,
$$
S(\mu_{(X_1,\dots,X_n)},\mu_{X_1}\otimes\dots\otimes\mu_{X_n})
:=\int\log{d\mu_{(X_1,\dots,X_n)}\over d(\mu_{X_1}\otimes\dots\otimes\mu_{X_n})}
\,d\mu_{(X_1,\dots,X_n)}
$$
if $\mu_{(X_1,\dots,X_n)}$ is absolutely continuous with respect to
$\mu_{X_1}\otimes\,\cdots\,\otimes\mu_{X_n}$; otherwise
$S(\mu_{(X_1,\dots,X_n)},\mu_{X_1}\otimes\dots\otimes\mu_{X_n}):=+\infty$.
When $H(X_i)>-\infty$ for all $1\le i\le n$, one can easily verify that
$$
S(\mu_{(X_1,\dots,X_n)},\mu_{X_1}\otimes\dots\otimes\mu_{X_n})
=-H(X_1,\dots,X_n)+\sum_{i=1}^nH(X_i).
$$
Thus, the above theorem yields the following:

\begin{cor}\label{C-1.7}
If $H(X_i)>-\infty$ for all $1\le i\le n$, then
\begin{align*}
I_\sym(X_1,\dots,X_n)&=\overline{I}_\sym(X_1,\dots,X_n) \\
&=S(\mu_{(X_1,\dots,X_n)},\mu_{X_1}\otimes\dots\otimes\mu_{X_n}).
\end{align*}
\end{cor}

\begin{cor}\label{C-1.8}
Under the same assumption as the above corollary, $I_\sym(X_1,\dots,X_n)=0$ if and
only if $X_1,\dots,X_n$ are independent.
\end{cor}

In particular, the original {\it mutual information} $I(X_1\wedge X_2)$ of two real random
variables $X_1,X_2$ is normally defined as
$$
I(X_1\wedge X_2):=S(\mu_{(X_1,X_2)},\mu_{X_1}\otimes\mu_{X_2}).
$$
Hence we have
$$
I(X_1\wedge X_2)=I_\sym(X_1,X_2)=\overline{I}_\sym(X_1,X_2)
$$
as long as $H(X_1)>-\infty$ and $H(X_2)>-\infty$ (and $X_1,X_2$ are bounded). For this
reason, we gave the term ``mutual information" to $I_\sym$.

Finally, some open problems are in order:
\begin{itemize}
\item[(1)] Without the assumption $H(X_i)>-\infty$ for $1\le i\le n$, does
$I_\sym(X_1,\dots,X_n)=\overline{I}_\sym(X_1,\dots,X_n)$ hold true?
\item[(2)] More strongly, does the limit such as
$$
\lim_{N\to\infty}{1\over N}\log\gamma_{S_N}^{\otimes n}
(\Delta_{\sym,R}(X_1,\dots,X_n;N,m,\delta))
$$
or
$$
\lim_{N\to\infty}{1\over N}\log\gamma_{S_N}^{\otimes n}
(\Delta_\sym(X_1,\dots,X_n:\xi_1(N),\dots,\xi_n(N);N,m,\delta))
$$
exist as in Lemma \ref{L-1.2}?
\item[(3)] Without the assumption $H(X_i)>-\infty$ for $1\le i\le n$, does
$I_\sym(X_1,\dots,X_n)
=S(\mu_{(X_1,\dots,X_n)},\mu_{X_1}\otimes\cdots\otimes\mu_{X_n})$ hold true?
Also, is $I_\sym(X_1,\dots,X_n)=0$ equivalent to the independence of $X_1,\dots,X_n$?
\item[(4)] Although the boundedness assumption for $X_1,\dots,X_n$ is rather essential
in the above discussions, it is desirable to extend the results in this section to
$X_1,\dots,X_n$ not necessarily bounded but having all moments.
\end{itemize}

\section{The discrete case}
\setcounter{equation}{0}

Let $\cY$ be a finite set with a probability measure $p$. The {\it Shannon entropy} 
of $p$ is
$$
S(p):=-\sum_{y \in\cY}p(y)\log p(y).
$$
For each sequence $\by=(y_1,\dots, y_N)\in\cY^N$, the {\it type} of $\by$ is a
probability measure on $\cY$ given by
$$
\nu_\by(t):=\frac{N_\by(t)}{N}\quad
\mbox{where}\quad N_\by(t):=\#\{j:y_j=t\},\quad t\in\cY.
$$
The number of possible types is smaller than $(N+1)^{\#\cY}$. 
If $\nu$ is a type and $\cT_N(\nu)$ denotes the set of all sequences of type $\nu$
from $\cY^N$, then the cardinality of $\cT_N(\nu)$ is estimated as follows:
\begin{equation}\label{F-2.1}
{1\over(N+1)^{\#\cY}}\,e^{NS(\nu)} \le \#\cT_N(\nu)\le e^{NS(\nu)}
\end{equation}
(see \cite[12.1.3]{CT} and \cite[Lemma 2.2]{CS}).

Let $p$ be a probability meausre on $\cY$. For each $N\in\bN$ and 
$\delta>0$ we define $\Delta(p;N,\delta)$ to be the set of all sequences 
$\by\in\cY^N$ such that $|\nu_\by(t)-p(t)|<\delta$ for all $t\in\cY$.
In other words, $\Delta(p;N,\delta)$ is the set of all $\delta$-typical
sequeces (with respect to the measure $p$). Then the next lemma is well known.

\begin{lemma}\label{L-2.1}
$$
S(p)=\lim_{\delta\searrow0}\lim_{N\to\infty}{1\over N}\log\#\Delta(p;N,\delta).
$$
\end{lemma}

In fact, this easily follows from \eqref{F-2.1}. Let $P_{N,\delta}$ be the
maximizer of the Shannon entropy on the set  of all types $\nu_\by$,
$\by\in\cY^N$, such that $|\nu_\by(t)-p(t)|<\delta$ for all $t\in\cY$. We can use the Shannon entropy of the type class corresponding to
$P_{N,\delta}$ to estimate the cardinality of $\Delta(p;N,\delta)$: 
$$
(N+1)^{-\#\cY} e^{NS(P_{N, \delta})}\le \#\Delta(p;N,\delta) 
\le e^{NS(P_{N,\delta})} (N+1)^{\#\cY}.
$$
It follows that 
\begin{align*}
\lim_{N\to\infty}{1\over N}\log\#\Delta(p;N,\delta)
&=\sup\{S(q):\mbox{$q$ is a probability meausre on $\cY$} \\
&\qquad\qquad\qquad\mbox{such that $|q(t)-p(t)|<\delta,\,t\in\cY$}\},
\end{align*}
and the lemma follows.

We consider the case where $p$ is the joint distribution of an $n$-tuple
$(X_1,\dots,X_n)$ of discrete random variables on $(\Omega,\bP)$. Throughout this
section we assume that the random variables $X_1,\dots,X_n$ have their values in
a finite set $\cX=\{t_1,\dots,t_d\}$.

\begin{definition}\label{D-2.2}{\rm
Let $p_{(X_1,\dots,X_n)}$ denote the joint distribution of $(X_1,\dots,X_n)$,
which is a measure on $\cX^n$ while the distribution $p_{X_i}$ of $X_i$ is a
measure on $\cX$, $1\le i\le n$. We write $\Delta(X_i;N,\delta)$ for
$\Delta(p_{X_i};N,\delta)$ and $\Delta(X_1,\dots,X_n;N,\delta)$ for
$\Delta(p_{(X_1,\dots,X_n)} ;N,\delta)$.
}\end{definition}

Next, we introduce the counterparts of Definitions \ref{D-1.3} and \ref{D-1.4} in
the discrete variable case.

\begin{definition}\label{D-2.3}{\rm
The action of $S_N$ on $\cX^N$ is similar to that on $\bR^N$ given in Defintion
\ref{D-1.3}. For $N\in\bN$ let $\cX_\le^N$ denote the set of all sequences of length
$N$ of the form
$$
\bx=(t_1,\dots,t_1,t_2,\dots,t_2,\dots,t_d,\dots,t_d).
$$
Oviously, such a sequence $\bx$ is uniquely determined by
$(N_\bx(t_1),\dots,N_\bx(t_d))$ or the type of $\bx$. That is, $\cX_\le^N$ is
regarded as the set of all types from $\cX^N$. For each $N\in\bN$ and $\delta>0$ we
denote by $\Delta_\sym(X_1,\dots,X_n;N,\delta)$ the set of all
$(\sigma_1,\dots,\sigma_n)\in S_N^n$ such that
$$
(\sigma_1(\bx_1),\dots,\sigma_n(\bx_n))\in\Delta(X_1,\dots,X_n;N,\delta)
$$
for some $(\bx_1,\dots,\bx_n)\in(\cX_\le^N)^n$. Define
$$
I_\sym(X_1,\dots,X_n):=-\lim_{\delta\searrow0}\limsup_{N\to\infty}
{1\over N}\log\gamma_{S_N}^{\otimes n}(\Delta_\sym(X_1,\dots,X_n;N,\delta)),
$$
and $\overline{I}_\sym(X_1,\dots,X_n)$ by replacing $\limsup$ by $\liminf$. Moreover,
for each $1\le i\le n$, choose a sequence $\xi_i=\{\xi_i(N)\}$ of
$\xi_i(N)=(\xi_i(N)_1,\dots,\xi_i(N)_N)\in\cX_\le^N$ such that
$\nu_{\xi_i(N)}\to p_{X_i}$ as $N\to\infty$. We then
define $\Delta_\sym(X_1,\dots,X_n:\xi_1(N),\dots,\xi_n(N);N,\delta)$,
$I_\sym(X_1,\dots,X_n:\xi_1,\dots,\xi_n)$ and
$\overline{I}_\sym(X_1,\dots,X_n:\xi_1,\dots,\xi_n)$ as in Definition \ref{D-1.4}.
}\end{definition}

\begin{lemma}\label{L-2.4}
For any choices of approximating sequences $\xi_1,\dots,\xi_n$ one has
\begin{align*}
I_\sym(X_1,\dots,X_n)&=I_\sym(X_1,\dots,X_n:\xi_1,\dots,\xi_n),\\
\overline{I}_\sym(X_1,\dots,X_n)
&=\overline{I}_\sym(X_1,\dots,X_n:\xi_1,\dots,\xi_n).
\end{align*}
\end{lemma}

\begin{proof}
It suffices to show that for each $\delta>0$ there are a $\delta'>0$ and an
$N_0\in\bN$ such that
\begin{equation}\label{F-2.2}
\Delta_\sym(X_1,\dots,X_n;N,\delta')
\subset\Delta_\sym(X_1,\dots,X_n:\xi_1(N),\dots,\xi_n(N);N,\delta)
\end{equation}
for all $N\ge N_0$. Choose $\delta'>0$ so that $3nd^{n+1}\delta'\le\delta$, where
$d=\#\cX$. Suppose $(\sigma_1,\dots,\sigma_n)$ is in the left-hand side of
\eqref{F-2.2} so that
$(\sigma_1(\bx_1),\dots,\sigma_n(\bx_n))\in\Delta(X_1,\dots,X_n;\allowbreak
N,\delta')$ for some $(\bx_1,\dots,\bx_n)$,
$\bx_i=(x_{i1},\dots,x_{iN})\in\cX_\le^N$. Since
\begin{equation}\label{F-2.3}
|\nu_{(\sigma_1(\bx_1),\dots,\sigma_n(\bx_n))}(z_1,\dots,z_n)
-p_{(X_1,\dots,X_n)}(z_1,\dots,z_n)|<\delta',
\quad(z_1,\dots,z_n)\in\cX^n,
\end{equation}
$$
\nu_{\bx_i}(t)=\sum_{z_1,\dots,z_{i-1},z_{i+1},\dots,z_n\in\cX}
\nu_{(\sigma_1(\bx_1),\dots,\sigma_n(\bx_n))}
(z_1,\dots,z_{i-1},t,z_{i+1},\dots,z_n),\quad t\in\cX,
$$
$$
p_{X_i}(t)=\sum_{z_1,\dots,z_{i-1},z_{i+1},\dots,z_n\in\cX}
p_{(X_1,\dots,X_n)}(z_1,\dots,z_{i-1},t,z_{i+1},\dots,z_n),\quad t\in\cX,
$$
it follows that
\begin{equation}\label{F-2.4}
|\nu_{\bx_i}(t)-p_{X_i}(t)|<d^{n-1}\delta'
\end{equation}
for any $1\le i\le n$ and $t\in\cX$. Now, choose an $N_0\in\bN$ so that
$|\nu_{\xi_i(N)}(t)-p_{X_i}(t)|<\delta'$ and hence
\begin{equation}\label{F-2.5}
|\nu_{\xi_i(N)}(t)-\nu_{\bx_i}(t)|<2d^{n-1}\delta'
\end{equation}
for any $1\le i\le n$ and $t\in\cX$ and for all $N\ge N_0$. Since
\begin{align*}
&|(N_{\xi_i(N)}(t_1)+\dots+N_{\xi_i(N)}(t_l))
-(N_{\bx_i}(t_1)+\dots+N_{\bx_i}(t_l))| \\
&\qquad\le|N_{\xi_i(N)}(t_1)-N_{\bx_i}(t_1)|
+\dots+|N_{\xi_i(N)}(t_l)-N_{\bx_i}(t_l)| \\
&\qquad<2Nd^n\delta'
\end{align*}
for every $1\le l\le d$ thanks to \eqref{F-2.5}, it is easily seen that
$$
\#\bigl\{j\in\{1,\dots,N\}:\xi_i(N)_j\ne x_{ij}\bigr\}<2Nd^{n+1}\delta'
$$
for any $1\le i\le n$. Hence we get
\begin{align*}
&|\nu_{(\sigma_1(\xi_1(N)),\dots,\sigma_n(\xi_n(N)))}(z_1,\dots,z_n)
-\nu_{(\sigma_1(\bx_1),\dots,\sigma_n(\bx_n))}(z_1,\dots,z_n)| \\
&\quad={1\over N}\big|\#\{j:\xi_1(N)_{\sigma_1^{-1}(j)}=z_1,\dots,
\xi_n(N)_{\sigma_n^{-1}(j)}=z_n\} \\
&\qquad\qquad-\#\{j:x_{1\sigma_1^{-1}(j)}=z_1,\dots,
x_{n\sigma_n^{-1}(j)}=z_n\}\big| \\
&\quad\le{1\over N}\sum_{i=1}^n\#\{j:\xi_i(N)_j\ne x_{ij}\}
<2nd^{n+1}\delta'
\end{align*}
so that thanks to \eqref{F-2.3}
$$
|\nu_{(\sigma_1(\xi_1(N)),\dots,\sigma_n(\xi_n(N)))}(z_1,\dots,z_n)
-p_{(X_1,\dots,X_n)}(z_1,\dots,z_n)|
<3nd^{n+1}\delta'\le\delta
$$
for every $(z_1,\dots,z_n)\in\cX^n$. Therefore, $(\sigma_1,\dots,\sigma_n)$ is in
the right-hand side of \eqref{F-2.2}, as required.
\end{proof}

The next theorem is the discrete variable version of Theorem \ref{T-1.6}.

\begin{thm}\label{T-2.5}
$$
I_\sym(X_1,\dots,X_n)=\overline{I}_\sym(X_1,\dots,X_n)
=-S(X_1,\dots,X_n)+\sum_{i=1}^nS(X_i).
$$
\end{thm}

\begin{proof}
For each sequence $(N_1,\dots,N_d)$ of integers $N_l\ge0$ with $\sum_{l=1}^dN_l=N$,
let $S(N_1,\dots,N_d)$ denote the subgroup of $S_N$ consisting of products of
permutations of $\{1,\dots,N_1\}$, $\{N_1+1,\dots,N_1+N_2\}$, $\dots$,
$\{N_1+\dots+N_{d-1}+1,\dots,N\}$, and let
$$
S_N/S(N_1,\dots,N_d)
$$
be the set of left cosets of $S(N_1,\dots,N_d)$. For each $\bx\in\cX_\le^N$ and
$\sigma\in S_N$ we write $[\sigma]_\bx$ for the left coset of
$S(N_\bx(t_1),\dots,N_\bx(t_d))$ containing $\sigma$. Then it is clear that every
$\bs\in\cX^N$ is represented as $\bs=\sigma(\bx)$ with a unique pair
$(\bx,[\sigma]_\bx)$ of $\bx\in\cX_\le^N$ and
$[\sigma]_\bx\in S_N/S(N_\bx(t_1),\dots,N_\bx(t_d))$.

For any $\eps>0$ one can choose a $\delta>0$ such that for every $1\le i\le n$ and
every probability measure $p$ on $\cX$, if $|p(t)-p_{X_i}(t)|<\delta$ for all
$t\in\cX$, then $|S(p)-S(p_{X_i})|<\eps$. This implies that for each $N\in\bN$ and
$1\le i\le n$, one has $|S(\nu_\bx)-S(p_{X_i})|<\eps$ whenever
$\bx\in\Delta(X_i;N,\delta)$. Notice that
$\Delta_\sym(X_1,\dots,X_n;N,\delta/d^{n-1})$ is the union of
$[\sigma_1]_{\bx_1}\times\dots\times[\sigma_n]_{\bx_n}$ for all
$(\bx_1,\dots,\bx_n;[\sigma_1]_{\bx_1},\dots,[\sigma_n]_{\bx_n})$ of
$\bx_i\in\cX_\le^N$ and
$[\sigma_i]_{\bx_i}\in S_N/S(N_{\bx_i}(t_1),\dots,N_{\bx_i}(t_d))$ such that
$(\sigma_1(\bx_1),\dots,\sigma_n(\bx_n))\in\Delta(X_1,\dots,X_n;N,\delta/d^{n-1})$.
Now, suppose $(\bx_1,\dots,\bx_n)\in(\cX_\le^N)^n$,
$(\sigma_1,\dots,\sigma_n)\in S_N^n$ and
$(\sigma_1(\bx_1),\dots,\sigma_n(\bx_n))\in\Delta(X_1,\dots,X_n;N,\delta/d^{n-1})$.
Then, for each $1\le i\le n$ we get $\bx_i\in\Delta(X_i;N,\delta)$, i.e.,
$|\nu_{\bx_i}(t)-p_{X_i}(t)|<\delta$ for all $t\in\cX$ as \eqref{F-2.4}.
Hence we have
\begin{equation}\label{F-2.6}
\#\bigl([\sigma_1]_{\bx_1}\times\dots\times[\sigma_n]_{\bx_n}\bigr)
\le\prod_{i=1}^n\Biggl(\max_{\bx\in\Delta(X_i;N,\delta)}
\prod_{t\in\cX}N_\bx(t)!\Biggr)
\end{equation}
so that
\begin{align*}
&\#\Delta_\sym(X_1,\dots,X_n;N,\delta/d^{n-1}) \\
&\qquad\le\#\Delta(X_1,\dots,X_n;N,\delta/d^{n-1})\cdot
\prod_{i=1}^n\Biggl(\max_{\bx\in\Delta(X_i;N,\delta)}
\prod_{t\in\cX}N_\bx(t)!\Biggr).
\end{align*}
Therefore,
\begin{align}\label{F-2.7}
&{1\over N}\log\gamma_{S_N}^{\otimes n}
\bigl(\Delta_\sym(X_1,\dots,X_n;N,\delta/d^{n-1})\bigr) \nonumber\\
&\qquad\le{1\over N}\log\#\Delta(X_1,\dots,X_n;N,\delta/d^{n-1}) \nonumber\\
&\qquad\quad+\sum_{i=1}^n\max_{\bx\in\Delta(X_i;N,\delta)}
\Biggl({1\over N}\sum_{t\in\cX}\log N_\bx(t)!\Biggr)
-{n\over N}\log N!.
\end{align}
For each $1\le i\le n$ and for any $\bx\in\Delta(X_i;N,,\delta)$, the Stirling
formula yields
\begin{align}\label{F-2.8}
&{1\over N}\sum_{t\in\cX}\log N_\bx(t)!-{1\over N}\log N! \nonumber\\
&\qquad=\sum_{t\in\cX}\biggl({N_\bx(t)\over N}\log N_\bx(t)
-{N_\bx(t)\over N}\biggr)-\log N+1+o(1) \nonumber\\
&\qquad=-S(\nu_\bx)+o(1)\le-S(p_{X_i})+\eps+o(1)
\quad\mbox{as $N\to\infty$}
\end{align}
thanks to the above choice of $\delta>0$. Here, note that the $o(1)$ in the above
estimate is uniform for $\bx\in\Delta(X_i;N,\delta)$.
Hence, by \eqref{F-2.7}, \eqref{F-2.8} and by Lemma \ref{L-2.1} applied to
$p_{(X_1,\dots,X_n)}$ on $\cX^n$, we obtain
$$
-I_\sym(X_1,\dots,X_n)\le S(p_{(X_1,\dots,X_n)})-\sum_{i=1}^nS(p_{X_i})+n\eps
$$
and hence
\begin{equation}\label{F-2.9}
I_\sym(X_1,\dots,X_n)\ge-S(X_1,\dots,X_n)+\sum_{i=1}^nS(X_i).
\end{equation}

Next, we prove the converse direction. For any $\eps>0$ choose a $\delta>0$ as above.
For $N\in\bN$ let $\Xi(N,\delta/d^{n-1})$ be the set of all
$(\bx_1,\dots,\bx_n)\in(\cX_\le^N)^n$ such that
$$
(\sigma_1(\bx_1),\dots,\sigma_n(\bx_n))\in\Delta(X_1,\dots,X_n;N,\delta/d^{n-1})
$$
for some $(\sigma_1,\dots,\sigma_n)\in S_N^n$. Furthermore, for each
$(\bx_1,\dots,\bx_n)\in\Xi(N,\delta/d^{n-1})$, let
$\Sigma(\bx_1,\dots,\bx_n;N,\delta/d^{n-1})$ be the set of all
$$
([\sigma_1]_{\bx_1},\dots,[\sigma_n]_{\bx_n})
\in\prod_{i=1}^nS_N/S(N_{\bx_i}(t_1),\dots,N_{\bx_i}(t_d))
$$
such that
$(\sigma_1(\bx_1),\dots,\sigma_n(\bx_n))\in\Delta(X_1,\dots,X_n;N,\delta/d^{n-1})$.
Then it is obvious that
\begin{equation}\label{F-2.10}
\#\Delta(X_1,\dots,X_n;N,\delta/d^{n-1})
\le\sum_{(\bx_1,\dots,\bx_n)\in\Xi(N,\delta/d^{n-1})}
\#\Sigma(\bx_1,\dots,\bx_n;N,\delta/d^{n-1}).
\end{equation}
When $(\bx_1,\dots,\bx_n)\in\Xi(N,\delta/d^{n-1})$, we get
$\bx_i\in\Delta(X_i;N,\delta)$ as \eqref{F-2.4} for $1\le i\le n$. Hence it is seen
that
\begin{align}\label{F-2.11}
\#\Xi(N,\delta/d^{n-1})
&\le\prod_{i=1}^n\#\Delta(X_i;N,\delta) \nonumber\\
&=\prod_{i=1}^n\#\bigl\{(N_1,\dots,N_d): N_l\ge0\ \mbox{is an integer in}
\nonumber\\
&\qquad\qquad\quad\bigl(N(p_{X_i}(t_l)-\delta),N(p_{X_i}(t_l)+\delta)\bigr)
\ \mbox{for $1\le l\le d$}\bigr\} \nonumber\\
&<(2N\delta+1)^{nd}.
\end{align}
For any fixed $(\bx_1,\dots,\bx_n)\in\Xi(N,\delta/d^{n-1})$, suppose
$([\sigma_1]_{\bx_1},\dots,[\sigma_n]_{\bx_n})\in\Sigma(\bx_1,\dots,\bx_n;\allowbreak
N,\delta/d^{n-1})$; then we get
$$
\#\bigl([\sigma_1]_{\bx_1}\times\dots\times[\sigma_n]_{\bx_n}\bigr)
\ge\prod_{i=1}^n\Biggl(\min_{\bx\in\Delta(X_i;N,\delta)}
\prod_{t\in\cX}N_\bx(t)!\Biggr)
$$
similarly to \eqref{F-2.6}. Therefore,
\begin{align}\label{F-2.12}
&\#\Delta_\sym(X_1,\dots,X_n;N,\delta/d^{n-1}) \nonumber\\
&\qquad\ge\sum_{([\sigma_1]_{\bx_1},\dots,[\sigma_n]_{\bx_n})
\in\Sigma(\bx_1,\dots,\bx_n;N,\delta/d^{n-1})}
\#\bigl([\sigma_1]_{\bx_1}\times\dots\times[\sigma_n]_{\bx_n}\bigr) \nonumber\\
&\qquad\ge\#\Sigma(\bx_1,\dots,\bx_n;N,\delta/d^{n-1})
\cdot\prod_{i=1}^n\Biggl(\min_{\bx\in\Delta(X_i;N,\delta)}
\prod_{t\in\cX}N_\bx(t)!\Biggr).
\end{align}
By \eqref{F-2.10}--\eqref{F-2.12} we obtain
$$
\#\Delta(X_1,\dots,X_n;N,\delta/d^{n-1})
\le{\#\Delta_\sym(X_1,\dots,X_n;N,\delta/d^{n-1})
\cdot(2N\delta+1)^{nd}\over
\prod_{i=1}^n\Bigl(\min_{\bx\in\Delta(X_i;N,\delta)}
\prod_{t\in\cX}N_\bx(t)!\Bigr)}
$$
so that
\begin{align*}
&{1\over N}\log\#\Delta(X_1,\dots,X_n;N,\delta/d^{n-1}) \\
&\qquad\le{1\over N}\log\gamma_{S_N}^{\otimes n}
\bigl(\Delta_\sym(X_1,\dots,X_n;N,\delta/d^{n-1})\bigr) \\
&\qquad\quad-\sum_{i=1}^n\min_{\bx\in\Delta(X_i;N,\delta)}\Biggl({1\over N}
\sum_{t\in\cX}\log N_\bx(t)!\Biggr)
+{n\over N}\log N!+{nd\over N}\log(2N\delta+1).
\end{align*}
Since it follows similarly to \eqref{F-2.8} that
$$
-{1\over N}\sum_{t\in\cX}\log N_\bx(t)!+{1\over N}\log N!
\le S(p_{X_i})+\eps+o(1)
\quad\mbox{as $N\to\infty$}
$$
with uniform $o(1)$ for all $\bx\in\Delta(X_i;N,\delta)$, we obtain
$$
S(p_{(X_1,\dots,X_n)})\le-\overline{I}_\sym(X_1,\dots,X_n)
+\sum_{i=1}^nS(p_{X_i})+n\eps
$$
by Lemma \ref{L-2.1} again, and hence
\begin{equation}\label{F-2.13}
\overline{I}_\sym(X_1,\dots,X_n)\le-S(X_1,\dots,X_n)
+\sum_{i=1}^n S(X_i).
\end{equation}
The conclusion follows from \eqref{F-2.9} and \eqref{F-2.13}.
\end{proof}

In particular, the mutual information $I(X_1\wedge X_2)$ of $X_1$ and $X_2$ is equivalently
expressed as
\begin{align*}
I(X_1\wedge X_2)&=S(p_{(X_1,X_2)},p_{X_1}\otimes p_{X_2})
=-S(p_{(X_1,X_2)})+S(p_{X_1})+S(p_{X_2}) \\
&=I_\sym(X_1,X_2)=\overline{I}_\sym(X_1,X_2).
\end{align*}

Similarly to the problem (2) mentioned in the last of Section 1, it is unknown
whether the limit
$$
\lim_{N\to\infty}{1\over N}\log
\gamma_{S_N}^{\otimes n}\bigl(\Delta_\sym(X_1,\dots,X_n;N,\delta)\bigr)
$$
exists or not.

\end{document}